\documentclass{asl}
\input{diagxy}

\title{On the equivalence of two quantifier elimination tests}

\author{Yimu Yin}
\revauthor{Yin, Yimu}

\address{Department of Philosophy\\
Carnegie Mellon University\\
Pittsburgh, PA 15217, USA}

\email{yimuy@andrew.cmu.edu}

\thanks{The author would like to thank J. Avigad, J. Cummings, and
R. Grossberg for many helpful discussions.}

\thanks{This work was partially supported by NSF grant DMS-0401042}

 \newtheorem{thm}{Theorem}[section]
 
 \newtheorem{lem}[thm]{Lemma}
 \newtheorem{fact}[thm]{Fact}

 \newtheorem{defn}[thm]{Definition}
 \newtheorem{ques}[thm]{Question}

 \numberwithin{equation}{section}

 \DeclareMathOperator{\cf}{cf}
 
 \DeclareMathOperator{\dom}{dom}

 \DeclareMathOperator{\id}{id}
 \DeclareMathOperator{\ed}{ED}
 \DeclareMathOperator{\cd}{CD}
 \DeclareMathOperator{\Th}{Th}
 \DeclareMathOperator{\tp}{tp}

 \newcommand{\abs}[1]{\left\vert#1\right\vert}
 \newcommand{\set}[1]{\left\{#1\right\}}
 \newcommand{\seq}[1]{\left<#1\right>}
 \newcommand{\norm}[1]{\left\Vert#1\right\Vert}

\newcommand{\limplies}{\rightarrow}

\newcommand{\liff}{\leftrightarrow}
\newcommand{\ex}[1]{\exists #1 \;} % exists x ...
\newcommand{\fa}[1]{\forall #1 \;} % forall x ...

\newcommand{\proves}{\vdash}

\newcommand{\fun}{\longrightarrow}

\newcommand{\sub}{\subseteq}
\newcommand{\esub}{\preceq}

\newcommand{\dash}{\mathalpha{\mbox{-}}} % e.g. for Sigma^1_1-AC

\newcommand{\T}{$T$\nobreakdash}

\newcommand{\aaa}{\bar{a}}
\newcommand{\bbb}{\bar{b}}
\newcommand{\ccc}{\bar{c}}
\newcommand{\ddd}{\bar{d}}

\newcommand{\x}{\bar{x}}
\newcommand{\y}{\bar{y}}

\newcommand{\tsub}{\mathbin{\mathchoice% truncated subtraction
{\buildrel .\lower.6ex\hbox{\vphantom{.}} \over {\smash-}}%
{\buildrel .\lower.6ex\hbox{\vphantom{.}} \over {\smash-}}%
{\buildrel .\lower.4ex\hbox{\vphantom{.}} \over {\smash-}}%
{\buildrel .\lower.3ex\hbox{\vphantom{.}} \over {\smash-}}}}

\begin{document}

\begin{abstract}
We prove that, for countable languages, two model-theoretic
quantifier elimination tests, one proposed by J.~R.~Shoenfield and
the other by L.~van~den~Dries, are equivalent.
\end{abstract}

\maketitle

%---------------------------------------------------------------------------------
\section{Introduction}
To facilitate the discussion we first introduce the following
terminological and notational conventions.

\begin{defn} Let $M$ be a model and $A \sub \abs{M}$. Let $N$ be the model $\seq{M, a}_{a \in A}$.
\begin{enumerate}
  \item The theory $\Th(N)$, denoted by $\cd(A, M)$, is called the \emph{complete
  diagram of $A$ in $M$}. If $A = \abs{M}$ we simply write $\cd(M)$.
  \item The set of all quantifier-free sentences in $\Th(N)$,
  denoted by $\ed(A, M)$, is called the \emph{elementary diagram of $A$ in
  $M$}. Again if $A = \abs{M}$ we simply write $\ed(M)$.
\end{enumerate}
\end{defn}

Obviously if $N \esub M$ then $\cd(N, M) = \cd(N)$ and if $N \sub M$
then $\ed(N, M) = \ed(N)$.

We say that a theory $T$ is \emph{model complete} if and only if,
for every pair of models $N, M \models T$, $N \sub M$ implies $N
\preceq M$. Abraham Robinson showed that under certain conditions a
model complete theory admits quantifier elimination (QE for short).
This was one of the results that inaugurated the use of
model-theoretic methods in the study of QE. Model-completeness has
many equivalent formulations:

\begin{fact}\label{subcom:modcom}
Let $T$ be any theory. The following are equivalent:
\begin{enumerate}
  \item\label{mc:mc} $T$ is model complete.
  \item\label{mc:two:models} For any two models $N, M \models T$ with $N \sub M$ there is an $N^* \models T$ such
  that $N \esub N^*$ and $M$ can be embedded into $N^*$ over $N$.
  \item\label{mc:complete} For any $M \models T$ the theory $T \cup \ed(M)$ is
  complete.
  \item\label{mc:existential:absolute} For any two models $N, M \models T$ with $N \sub M$, every existential formula $\varphi(\bar{x})$,
  and every $\bar{b} \in \abs{N}$, we have $M \models \varphi(\bar{b})$ if and only if $N \models
  \varphi(\bar{b})$.
  \item\label{mc:ex:uni} For every existential formula $\varphi(\bar{x})$ there is a
  universal formula $\varphi^*(\bar{x})$ such that $T \proves \varphi(\bar{x}) \liff
  \varphi^*(\bar{x})$.
  \item\label{mc:all:uni} For every formula $\varphi(\bar{x})$ there is a universal formula $\varphi^*(\bar{x})$
  such that $T \proves \varphi(\bar{x}) \liff \varphi^*(\bar{x})$.
  \item\label{mc:all:uni:ex} For every formula $\varphi(\bar{x})$ there is a universal formula
  $\varphi_1(\bar{x})$ and an existential formula $\varphi_2(\bar{x})$
  such that $T \proves \varphi_1(\bar{x}) \liff \varphi(\bar{x}) \liff \varphi_2(\bar{x})$.
\end{enumerate}
\end{fact}

For a proof of this fact see \cite{chang:keisler:90} and
\cite{Sac72}.

However, there are theories which are model complete but do not
admit QE. For example, the complete theory of real closed fields in
the language of rings is model complete, but the formula $\ex{x} x
\times x = y$ is not equivalent to any quantifier-free formula in
this theory. See~\cite{chang:keisler:90} for details.

Over the years many model-theoretic properties have been proposed to
strengthen model-completeness so that QE is implied without any
additional assumptions on the theory in question. Some of these
properties are logically equivalent to QE; others are strictly
stronger than QE. Below we shall prove that two of the stronger
ones, one proposed by J.~R.~Shoenfield and the other by
L.~van~den~Dries, are equivalent for countable languages.

%---------------------------------------------------------------------------------
\section{Some QE tests}
Let $T$ be any theory. Here are some model-theoretic QE tests that
are stronger than model-completeness:

\begin{defn}\label{QEtests:sc}
$T$ is \emph{submodel complete} if and only if for any model $M
\models T$ and any $N \subseteq M$ the theory $T \cup \ed(N)$ is
complete.
\end{defn}

This is a direct strengthening of
\ref{subcom:modcom}.\ref{mc:complete}.

\begin{defn}\label{QEtests:sa}
$T$ has the \emph{submodel amalgamation property} (SA-property for
short) if and only if for any $M_1, M_2 \models T$ and any $N \sub
M_1, M_2$ there is an $M^* \models T$ such that $M_1 \esub M^*$ and
$M_2$ can be embedded into $M^*$ over $N$ via a monomorphism $f$;
that is, the following diagram
\[
  \bfig
  \square/->`<-`<-`->/<450,450>[M_1`M^*`N`M_2;\esub`\sub`f`\sub]
  \efig
\]
commutes.
\end{defn}

This is a direct strengthening of
\ref{subcom:modcom}.\ref{mc:two:models}.

\begin{defn}\label{QEtests:s}
$T$ has the \emph{Shoenfield property} (S-property for short) if and
only if for any two models $M_1, M_2 \models T$ such that $M_2$ is
$\norm{M_1}^+$-saturated and any isomorphism $f: N_1 \fun N_2$ with
$N_1 \sub M_1$ and $N_2 \sub M_2$, there is a monomorphism $f^*: M_1
\fun M_2$ extending $f$.
\end{defn}

\begin{defn}\label{QEtests:ss}
$T$ has the \emph{strong Shoenfield property} (SS-property for
short) if and only if
    \begin{enumerate}
      \item For every two models $M_1, M_2 \models T$ and every
      two models $N_1 \sub M_1$ and $N_2 \sub M_2$, if $f: N_1 \fun N_2$ is an isomorphism,
      then there is an isomorphism $f^*: N_1^* \fun N_2^*$ which is an
      extension of $f$, where $N_1^* \subseteq M_1$, $N_2^* \subseteq M_2$, and $N_1^*, N_2^* \models T$;
      \item For every two models $N, M \models T$ with $N \subseteq M$, every existential formula $\varphi(\bar{x})$, and
      every $\bar{b} \in \abs{N}$, we have $M \models \varphi(\bar{b})$ if
      and only if $N \models \varphi(\bar{b})$. In other words, $T$
      is model complete.
    \end{enumerate}
\end{defn}

When there is no danger of confusion we abuse $L(T)$ to denote both
the language of $T$ and the set of all well-formed formulas in the
language of $T$. For two structures $N$ and $M$ in $L(T)$ we say
that $M$ is a \emph{$T$-extension} of $N$ if $\abs{N} \sub \abs{M}$
and $M \models T$.

\begin{defn}\label{QEtests:d}
$T$ has the \emph{van den Dries property} (D-property for short) if
and only if
    \begin{enumerate}
      \item For any model $N$, if there exists a model $M \models T$ such that
      $N \sub M$, then there is a $T$-closure $N^*$ of $N$, that is, a model
      $N^* \models T$ such that $N \sub N^*$ and $N^*$ can be embedded over $N$ into any $T$-extension
      of $N$;
      \item\label{dproperty:2} If $N, M \models T$ and $N \subsetneq M$, then
      there is an $a \in \abs{M} \setminus \abs{N}$ such that $N + a$
      can be embedded into an elementary extension of $N$ over $N$, where $N + a$
      is the smallest submodel of $M$ that contains $\abs{N} \cup
      \{a\}$.
    \end{enumerate}
\end{defn}

The SS-property first appeared in Shoenfield's
textbook~\cite{Shoen67}. He subsequently modified it into the
S-property and proved its equivalence to QE in~\cite{Shoen71}. The
D-property was given by van den Dries in~\cite{Dri85}
and~\cite{Dri88}, which is a straightforward strengthening of the
SS-property. However, the main result Theorem~\ref{D:SS:equi} below
shows that, for countable languages, its main advantage over the
SS-property is its conceptual concreteness rather than its logical
strength.

\begin{thm}\label{QEtests:relations}
Let $T$ be a theory in a language with at least one constant symbol.
For the following statements,
\begin{enumerate}
  \item\label{qe:sc} $T$ is submodel complete,
  \item\label{qe:sa} $T$ has the SA-property,
  \item\label{qe:sp} $T$ has the S-property,
  \item\label{qe:wsp} $T$ has the SS-property,
  \item\label{qe:dp} $T$ has the D-property,
  \item\label{qe:qe} $T$ admits QE,
\end{enumerate}
these logical implications hold:
\[
  \bfig
  \Ctriangle/<=``<=/<180,180>[\ref{qe:qe}`\ref{qe:sc}`\ref{qe:sa};``]
  \Dtriangle(180,0)/`=>`=>/<180,180>[\ref{qe:qe}`\ref{qe:sp}`\ref{qe:sa};``]
  \morphism(360,180)/<=/<210,0>[\ref{qe:sp}`\ref{qe:wsp};]
  \morphism(570,180)/<=/<210,0>[\ref{qe:wsp}`\ref{qe:dp};]
  \efig
\]
\end{thm}

\begin{proof}
That \ref{qe:sc}, \ref{qe:sa}, and \ref{qe:sp} are equivalent to QE
is well-known. See, for example, \cite{Sac72} and \cite{Shoen71}.
Here we give proofs to the remaining two implications. We also show
directly how the first condition of the SS-property achieves QE on
top of model-completeness. This proof is a modification of the
standard proof of
``\ref{subcom:modcom}.\ref{mc:existential:absolute}~$\Rightarrow$~\ref{subcom:modcom}.\ref{mc:ex:uni}''
in the literature, which establishes a crucial connection between
model-theoretic properties and syntactical properties.

\ref{qe:wsp} $\Rightarrow$~\ref{qe:qe}: Let $\varphi(\x)$ be a
formula in $L(T)$. Since $T$ is model complete,
by~\ref{subcom:modcom}, $\varphi(\x)$ is equivalent to both a
universal formula and an existential formula. Hence we may assume
that $\varphi(\x)$ is a universal formula. Let $\varphi^*(\x)$ be an
existential formula such that $T \proves \varphi(\x) \liff
\varphi^*(\x)$. Let $\ccc$ be new constants. Let $\Gamma$ be a set
that contains exactly the following formulas:
\begin{itemize}
  \item $T \cup \{\varphi(\ccc)\}$, and
  \item every quantifier-free formula $\neg \psi(\ccc)$ such that $T \proves \fa{\x} (\psi(\x) \limplies
  \varphi(\x))$.
\end{itemize}

Suppose for contradiction that $\Gamma$ is consistent. Take any
model $M \models \Gamma$. Let $N \sub M$ be the minimal submodel
generated by $\ccc$. Note that every element in $N$ can be written
as a term that only involves $\ccc$, the constants of $L(T)$, and
the functions of $L(T)$. Now, if $T \cup \ed(N)$ does not prove
$\varphi(\ccc)$, then fix a model $M^* \models T \cup \ed(N) \cup
\{\neg \varphi(\ccc)\}$. By the first condition of the SS-property
we can find an $N_1 \models T \cup \ed(N)$ in $M$ and an $N_2
\models T \cup \ed(N)$ in $M^*$ such that they are isomorphic over
$N$. Since $\varphi(\x)$ is a universal formula and $M \models
\varphi(\ccc)$, we have $N_1 \models \varphi(\ccc)$. So $N_2 \models
\varphi(\ccc)$, so $N_2 \models \varphi^*(\ccc)$, so $M^* \models
\varphi^*(\ccc)$, so $M^* \models \varphi(\ccc)$, contradiction. So
$T \cup \ed(N) \proves \varphi(\ccc)$. So there is a quantifier-free
formula $\psi(\ccc) \in \ed(N)$ such that $T \cup \{\psi(\ccc)\}
\proves \varphi(\ccc)$, so $T \proves \psi(\ccc) \limplies
\varphi(\ccc)$. But $\ccc$ are new constants, so $T \proves \fa{\x}
(\psi(\x) \limplies \varphi(\x))$. So $\neg \psi(\ccc) \in \Gamma$,
contradiction again.

So $\Gamma$ is not consistent. This means that there are finitely
many quantifier-free formulas $\psi_i(\x)$ such that $T \proves
\fa{\x} (\psi_i(\x) \limplies \varphi(\x))$ for every $i$ and $T
\proves \fa{\x} (\varphi(\x) \limplies \bigvee_i \psi_i(\x))$. So $T
\proves \fa{\x} (\varphi(\x) \liff \bigvee_i \psi_i(\x))$, as
desired.

\ref{qe:wsp} $\Rightarrow$~\ref{qe:sp}: Let $M_1, M_2 \models T$, $N
\sub M_1, M_2$, and let $M_2$ be $\norm{M_1}^+$-saturated. By the
first condition of the SS-property we can find two \T-extensions
$N_1, N_2$ of $N$ in $M_1, M_2$ respectively that are isomorphic
over $N$. Let the isomorphism be $f$. Pick an $a \in \abs{M_1}
\setminus \abs{N_1}$ and consider any quantifier-free formula
$\varphi(x; \bbb)$ with $\bbb \in \abs{N_1}$ such that $M_1 \models
\varphi(a; \bbb)$. Since $M_1 \models \ex{x} \varphi(x; \bbb)$, by
the second condition of the SS-property we have $N_1 \models \ex{x}
\varphi(x; \bbb)$, so $N_2 \models \ex{x} \varphi(x; f(\bbb))$, so
$M_2 \models \ex{x} \varphi(x; f(\bbb))$. Hence the quantifier-free
type $f(p)$ is realized in $M_2$, say, by $d$, where $p$ is the set
of all quantifier-free formulas in $\tp(a / \abs{N_1}, M_1)$. If we
set $a \longmapsto d$ then we get an induced isomorphism between
$N_1 + a$ and $N_2 + d$. Iterating this procedure to exhaust all
elements in $M_1$ we see that $M_1$ can be embedded into $M_2$ over
$N$.

\ref{qe:dp} $\Rightarrow$~\ref{qe:wsp}: Trivially the closure
property, that is, the first condition of the D-property, implies
the first condition of the SS-property. For the second condition of
the SS-property, let $N, M \models T$ with $N \subseteq M$. Consider
an existential formula $\ex{\x} \varphi(\x; \bbb)$ that is satisfied
in $M$, where $\bar{b} \in \abs{N}$ and $\varphi(\x; \bbb)$ is
quantifier-free. So let $\ccc$ be such that $M \models \varphi(\ccc;
\bbb)$. We construct the following diagram:
\[
  \bfig
  \Vtrianglepair<450,450>[N_0`N_0 + a_0`N_1`N_0^*;\sub`\sub`\preceq``f_0]
  \Vtrianglepair(900,0)<450,450>[N_1`N_1 + a_1`N_2`N_1^*;\sub`\sub`\preceq``f_1]
  \morphism(1800,450)<450,-450>[N_2`\cdots;]
  \morphism(1800,450)/.>/<900,0>[N_2`M;\sub]
  \efig
\]
where $N_0 = N$, each $N_{i+1}$ is the \T-closure of $N_i + a_i$
promised by the closure property, each $a_i$ and $N_i^*$ are as
described in the second condition of the D-property, all arrows are
monomorphisms, and at the limit stage we simply take the union of
all previous $N_i$'s.

Now, let $i$ be the least index such that $\ccc \in N_i$. Note that
$i$ cannot be a limit ordinal. So $N_i \models \ex{\x} \varphi(\x;
\bbb)$, so $N_{i-1}^* \models \ex{\x} \varphi(\x; \bbb)$, so
$N_{i-1} \models \ex{\x} \varphi(\x; \bbb)$, etc. If $\gamma$ is a
limit ordinal and $N_{\gamma} \models \ex{\x} \varphi(\x; \bbb)$,
then there is a $\ddd \in \abs{N_{\gamma}}$ such that $N_{\gamma}
\models \varphi(\ddd; \bbb)$, so by the construction there is a $j <
\gamma$ such that $\ddd \in \abs{N_j}$, so $N_j \models
\varphi(\ddd; \bbb)$, so $N_j \models \ex{\x} \varphi(\x; \bbb)$. As
we trace back in the diagram we see that $N = N_0 \models \ex{\x}
\varphi(\x; \bbb)$.
\end{proof}

The reason that we have assumed that the language of $T$ has at
least one constant symbol is to avoid certain pathology. That is, in
the proof of ``\ref{qe:wsp} $\Rightarrow$~\ref{qe:qe}'' above, if
$\varphi$ is a sentence and $L(T)$ has no constant symbol, then
$\ccc$ is the empty sequence and cannot generate any submodel as we
do not allow an empty model. The reader should observe that in this
case the proof will not go through if we simply use an arbitrary
submodel. In the sequel we shall always assume that $T$ has a
constant symbol whenever we are in a similar situation.

There are still more model-theoretic tests that are equivalent to
QE. They are all more or less variations of the three equivalent
tests in the above theorem. See~\cite{Hodges93} for more details
about this. On the other hand, it is tempting to ask if in the above
theorem all of the statements are indeed equivalent.

Jeremy Avigad has an example which shows that QE is strictly weaker
than the SS-property. Consider the set $2^{\omega}$ of all binary
sequences of length $\omega$. For each $n \in \omega$ let $Z_n$ be a
unary predicate such that if $n = 0$ then $Z_n(\eta)$ for any $\eta
\in 2^{\omega}$, otherwise $Z_n(\eta)$ if and only if $(\eta)_n =
0$. Let $T = \Th(\seq{2^{\omega}, Z_n}_{n \in \omega})$. Since
except equality all predicates in the language are unary, every
existential formula $\ex{x} \varphi(x; \y)$ is equivalent to a
formula of the form $\bigvee_i (\theta_i(\y) \wedge \ex{x} \phi_i(x;
\y))$, where $\phi_i(x; \y)$ is a conjunction of literals each of
which contains $x$. If the unary predicates in the formula $\ex{x}
\phi_i(x; \y)$ describe a ``consistent'' finite sequence, then it
can be translated into an equivalent quantifier-free formula that
only involves $\y$. So $T$ proves that every existential formula is
equivalent to a quantifier-free formula, which means that $T$ admits
QE. Now, it is not hard to see that any dense subset of $2^{\omega}$
is a model of $T$. Let $S_0 \sub 2^{\omega}$ be the set of those
sequences that have only finitely many 0's. Let $S_1 \sub
2^{\omega}$ be the set of those sequences that have only finitely
many 1's and the constant sequence $\bar{1}$. So both $S_0$ and
$S_1$ are models of $T$. Notice that $\{\bar{1}\}$ is a submodel of
both models as there is no function symbol in the language. Clearly
there cannot be isomorphic \T-extensions of $\{\bar{1}\}$ in $S_0$
and $S_1$.

What about the SS-property and the D-property? First of all it is
trivial that if a theory $T$ admits QE then the second condition of
the D-property holds, because, by \ref{subcom:modcom}, if $N, M
\models T$ and $N \sub M$ then $M$ itself is an elementary extension
of $N$. The closure property, however, is much harder to achieve.
The rest of this paper is devoted to proving

\begin{thm}\label{D:SS:equi}
For countable languages the SS-property and the D-property are
equivalent.
\end{thm}

The argument is by a transfinite induction.

%---------------------------------------------------------------------------------
\section{The base case of the induction}
We need more concepts and Henkin's Omitting Type Theorem.

\begin{defn}
Let $\bar{x}$ be a sequence of variables and $p$ a \T-type in
$\bar{x}$. If there exists a formula $\varphi(\x)$ such that $T \cup
\{\varphi(\x)\}$ is consistent and $\varphi(\x) \proves p$, then we
say that $p$ is \emph{isolated by $\varphi(\x)$ via $T$}. If in
context it is clear that which theory is being discussed then we
omit $T$.
\end{defn}

Note that if $p$ is a complete \T-type then $p$ is isolated via $T$
if and only if there exists a $\varphi \in p$ such that $\varphi
\proves p$.

\begin{defn}
Let $M \models T$ and $A \sub \abs{M}$. We say that $M$ is
\emph{almost \T-primary over $A$} if there exists an ordinal
$\alpha$ and a sequence $\seq{(N_i, b_i): i < \alpha}$ such that
\begin{enumerate}
  \item $N_0$ is the minimal submodel of $M$ that contains $A$,
  \item $b_i \in \abs{M} \setminus \abs{N_i}$ and $N_{i+1} = N_i + b_i$ for each $i < \alpha$ (if $\alpha = \beta + 1$
  then $b_{\beta}$ is not defined),
  \item $N_{\beta} = \bigcup_{i < \beta} N_i$ if $\beta$ is a limit ordinal and
  $\bigcup_{i < \alpha} N_i = M$,
  \item the type $\tp(b_j/\abs{N_j}, M)$ is isolated via $T_j$ for every $j < \alpha$,
  where $T_j = T \cup \cd(N_j, M)$.
\end{enumerate}
The sequence $\seq{(N_i, b_i): i < \alpha}$ is called an
\emph{almost isolating sequence for $M$ over $A$}. The ordinal
$\alpha$ is the \emph{length} of the sequence.
\end{defn}

For convenience, if $T = \Th(M)$ then we omit $T$. Also, sometimes
we allow an almost isolating sequence to have repeated consecutive
$b_i$'s. Of course in this case we no longer require $b_i \notin
\abs{N_i}$ for the repeated occurrences. Note that this definition
is a variation of the notion of a primary model, which plays an
important role in the proof of Morley's Theorem.

\begin{defn}
Let $M \models T$ and $A \sub \abs{M}$. We say that $M$ is
\emph{\T-primary over $A$} if there exists an ordinal $\alpha$ and
an enumeration $\seq{b_i: i < \alpha}$ of $\abs{M} \setminus A$ such
that the type
\[
\tp(b_j/A \cup \set{b_i: i < j}, M)
\]
is isolated via $T_j$ for every $j < \alpha$, where $T_j = T \cup
\cd(A \cup \set{b_i: i < j}, M)$. The sequence $\seq{b_i: i <
\alpha}$ is called an \emph{isolating sequence for $M$ over $A$}.
The ordinal $\alpha$ is the \emph{length} of the sequence.
\end{defn}

It is not hard to see that if $T$ is submodel complete and $N \sub M
\models T$ then $M$ is almost \T-primary over $N$ if and only if $M$
is \T-primary over $N$. We prefer the concept of an almost primary
model below because it is more explicit about what property is being
exploited, namely submodel completeness.

\begin{thm}[Henkin's Omitting Type Theorem]\label{henkin:omitting}
If $L(T)$ is countable and $\Gamma$ is a countable collection of
$T$-types such that $p$ is not isolated for every $p \in \Gamma$,
then there exists a countable model $M \models T$ that omits all the
types in $\Gamma$.
\end{thm}

We proceed to develop a couple of technical lemmas. We have the
following basic fact about an almost primary model satisfying a
submodel complete theory:

\begin{lem}\label{primary:basic}
Suppose $T$ is submodel complete. Let $N \sub M \models T$. Then: if
$M$ is almost \T-primary over $N$, then for every model $M^* \models
T \cup \ed(N)$ there is an elementary embedding from $M$ into $M^*$
over $N$.
\end{lem}

\begin{proof}
Since $T$ is submodel complete, the theory $T \cup \ed(N)$ is
complete. This means that for any formula $\varphi(\x)$ and any
$\aaa \in \abs{N}$ we have
\[
M \models \varphi(\aaa) \text{ iff } M^* \models \varphi(\aaa).
\]
Let $\seq{(N_i, b_i): i < \alpha}$ be an almost isolating sequence
for $M$ over $N$. So by definition $N_0 = N$. In order to prove the
lemma it is enough to construct a continuous sequence of
monomorphisms $g_i: N_i \fun M^*$ for $i < \alpha$ such that
\begin{enumerate}
  \item $g_0 = \id_{N}$,
  \item $N_i \models \varphi(\aaa)$ iff $M^* \models
  \varphi(g_i(\aaa))$ for each formula $\varphi(\x)$ and each
  $\aaa \in N_i$,
  \item if $i < j < \alpha$ then $g_i \sub g_j$, and
  \item if $\beta$ is a limit then $g_{\beta} = \bigcup_{i < \beta}
  g_i$.
\end{enumerate}
The embedding $g = \bigcup_{i < \alpha} g_i$ is as desired. That $g$
is elementary is because submodel completeness implies model
completeness (see \ref{subcom:modcom} and \ref{QEtests:relations}).

Now we proceed to construct the sequence. Due to the clause 4 all we
have to do is to make the successor case work. So suppose we have
successfully constructed the sequence up to the ordinal $i <
\alpha$. Since the complete type $p_i = \tp(b_i/\abs{N_i}, M)$ is
isolated via $T_i$ where $T_i = T \cup \cd(N_i, M)$, there exists a
formula $\varphi(x; \aaa) \in p_i$ isolating it. By the clause 2 we
have
\begin{equation}
\varphi(x; \aaa) \proves p_i \Rightarrow \varphi(x; g_i(\aaa))
\proves g_i(p_i). \tag{$\star$}\label{primary:type:transfer}
\end{equation}
Since $M \models \varphi(b_i; \aaa)$, we have $M \models \ex{x}
\varphi(x; \aaa)$, so $M^* \models \ex{x} \varphi(x; g_i(\aaa))$.
Let $c_i \in \abs{M^*}$ such that $M^* \models \varphi(c_i;
g_i(\aaa))$. So by~\eqref{primary:type:transfer} $c_i$ realizes the
type $g_i(p_i)$. Now define a function $g_{i + 1}$ by setting
$\tau(b_i) \longmapsto \tau(c_i)$ for each term $\tau(x)$ of
$L(T_i)$. It is easy to see that this is a well-defined monomorphism
from $N_{i+1}$ into $M^*$ which extends $g_i$ and takes $b_i$ to
$c_i$. That the clause 2 is satisfied is, again, because $T$ is
submodel complete.
\end{proof}

In order to build almost primary models we need the next crucial
lemma.

\begin{lem}\label{countable:submodel:isolated}
Suppose that $L(T)$ is countable and $T$ has the SS-property. Then
for
\begin{enumerate}
  \item every model $M \models T$,
  \item every countable submodel $N \sub M$,
  \item every formula $\varphi(x; \y)$ and every $\aaa \in \abs{N}$ such that
  $\ex{x} \varphi(x; \aaa) \in T \cup \ed(N)$ but $M \models \neg \varphi(b; \aaa)$ for every $b \in
  \abs{N}$,
\end{enumerate}
there is an element $c \in \abs{M} \setminus \abs{N}$ such that the
type $\tp(c/\abs{N}, M)$ is isolated and $M \models \varphi(c;
\aaa)$.
\end{lem}

\begin{proof}
Fix an $M$, an $N$, an $\aaa$, and a $\varphi(x; \y)$ as above.
Without loss of generality we may assume $M$ is countable as well.
Since $T$ has the SS-property, by~\ref{QEtests:relations}, the
theory $T \cup \ed(N)$ is complete. So $M \models \ex{x} \varphi(x;
\aaa)$. So $\varphi(M; \aaa) \neq \emptyset$ and, by the third
condition, $\varphi(M; \aaa) \sub \abs{M} \setminus \abs{N}$, where
$\varphi(M; \aaa)$ is the set $\set{c \in \abs{M}: M \models
\varphi(c; \aaa)}$. Also note that $T$ is model complete.

Suppose for contradiction we cannot find an element $c$ in $M$ as
required. Define a collection $\Gamma$ of $T \cup \ed(N)$-types:
\[
\Gamma = \set{\tp(c/\abs{N}, M): c \in \abs{M} \setminus \abs{N}
\text{ and } M \models \varphi(c; \aaa)}.
\]
Since $\Gamma$ is countable, by Henkin's Omitting Type Theorem there
is a model $O \models T \cup \ed(N)$ that omits every type in
$\Gamma$. But $T$ has the SS-property, so we can find two models
$M^* \sub M$, $O^* \sub O$ of $T$ such that there is an isomorphism
$h: M^* \cong O^*$ whose restriction to $N$ is $\id_N$. Since
$\ex{x} \varphi(x; \aaa) \in T \cup \ed(N)$, there must be some $c
\in \abs{M^*} \setminus \abs{N}$ such that $M^* \models \varphi(c;
\aaa)$. Since $T$ is model complete, we deduce
\[
\varphi(x; \aaa) \in \tp(c/\abs{N}, M^*) = \tp(c/\abs{N}, M).
\]
This means that $h(c)$ realizes the $T \cup \ed(N)$-type
$\tp(c/\abs{N}, M)$ in $O$, contradicting the choice of $O$.
\end{proof}

Note that in the above lemma, if $N$ is not a model of $T$, then
there must exist a formula $\ex{x} \varphi(x; \aaa) \in T \cup
\ed(N)$ with $\aaa \in \abs{N}$ such that $M \models \neg \varphi(b;
\aaa)$ for every $b \in \abs{N}$, because otherwise $N$ would be a
model of $T$ by the Tarski-Vaught Test as $T \cup \ed(N)$ is
complete. This property is important for our argument. We shall give
it a name:

\begin{defn}
Let $M \models T$, $N \sub M$, and $\aaa \in \abs{N}$. We say that
$\varphi(x; \aaa)$ is \emph{critical} for $N$ if $\ex{x} \varphi(x;
\aaa) \in T \cup \ed(N)$ and $\varphi(M; \aaa) \sub \abs{M}
\setminus \abs{N}$.
\end{defn}

Now the SS-property enables us to construct almost primary models
over countable submodels.

\begin{thm}\label{ss:countable:closure}
If $L(T)$ is countable and $T$ has the SS-property then, for any
model $M \models T$ and any countable submodel $N \sub M$, $N$ has a
$T$-closure.
\end{thm}

\begin{proof}
Fix $N \sub M \models T$ such that $N$ is countable. Again we may
assume that $M$ is countable as well. So by
Lemma~\ref{primary:basic} all we need to do is to build an almost
\T-primary model $N^*$ over $N$ inside $M$. For this it is enough to
build an almost isolating sequence for some model of $T$ over $N$.
The idea here is of course to find a suitable Skolem hull of $N$
inside $M$ such that the type of each ``key'' new element we find is
isolated over all the previous elements.

To be precise, we want to build an almost isolating sequence
$\seq{(N_i, b_i): i < \omega \cdot \omega}$ over $N$ such that for
\begin{itemize}
  \item each $n < \omega$,
  \item each $\aaa \in N_{\omega \cdot n}$, and
  \item each formula $\varphi(x; \y)$ such that $M \models \ex{x} \varphi(x;
  \aaa)$,
\end{itemize}
there is an $m < \omega$ such that $M \models \varphi(\tau(b_{\omega
\cdot n + m}); \aaa)$ for some term $\tau(x)$ in the language $L(T
\cup \ed(N_{\omega \cdot n + m}))$. It should be clear that
$\bigcup_{i < \omega \cdot \omega} N_i = N^*$ is an elementary
submodel of $M$, and hence is almost \T-primary over $N$.

Now we carry out the construction. Start with $N_0 = N$ of course.
Suppose $\seq{(N_i, b_i): i < \omega \cdot n}$ is defined. Let
$\seq{\varphi_k(x; \aaa_k): k < \omega}$ be an enumeration of all
the formulas in $T \cup \ed(N_{\omega \cdot n})$ such that for every
$k < \omega$ we have $M \models \ex{x} \varphi_k(x; \aaa_k)$ but $M
\models \neg \varphi_k(d; \aaa_k)$ for every $d \in N_{\omega \cdot
n}$. Now suppose we have extended the sequence all the way up to
$(N_{\omega \cdot n + k}, b_{\omega \cdot n + k})$ for some $k <
\omega$. Let $N_{\omega \cdot n + k + 1} = N_{\omega \cdot n + k} +
b_{\omega \cdot n + k}$. If there is a $d \in N_{\omega \cdot n + k
+ 1}$ such that $M \models \varphi_{k+1}(d; \aaa_{k+1})$ then let
$b_{\omega \cdot n + k + 1} = b_{\omega \cdot n + k}$. Otherwise by
Lemma~\ref{countable:submodel:isolated} we can pick a $b_{\omega
\cdot n + k + 1} \in \abs{M} \setminus \abs{N_{\omega \cdot n + k +
1}}$ such that $M \models \varphi_{k+1}(b_{\omega \cdot n + k + 1};
\aaa_{k+1})$ and the type $\tp(b_{\omega \cdot n + k +
1}/\abs{N_{\omega \cdot n + k + 1}}, M)$ is isolated.
\end{proof}

%---------------------------------------------------------------------------------
\section{The inductive step}
The reader may ask: What is preventing us here from simply extending
the above theorem to arbitrary theories and arbitrary submodels? One
difficulty is this: We do not know how to extend Henkin's Omitting
Type Theorem to uncountable languages and hence are unable to
develop an analog of Lemma~\ref{countable:submodel:isolated} for
uncountable languages. In fact if we simply drop the countability
requirement in Henkin's Omitting Type Theorem then it is false.
See~\cite{chang:keisler:90} for discussions. However, in this last
section we will show how to circumvent this difficulty if the
language in question is countable. For this we need some basic
concepts and facts in infinitary combinatorics, in particular
stationary sets and Fodor's Lemma.

Throughout the rest of this section $T$ is a theory in a countable
language and has the SS-property. Our strategy is to establish an
analog of Lemma~\ref{countable:submodel:isolated} for any submodel.
Let $M \models T$ and $N \sub M$ such that $N$ is uncountable and is
not a model of $T$. We have two cases to consider, namely $\norm{N}$
is regular and $\norm{N}$ is singular.

\begin{defn}
Let $\alpha$ be an ordinal. A sequence $\seq{N_i: i < \alpha}$ is an
\emph{$\alpha$-resolution of $N$} if
\begin{enumerate}
  \item $N_i$ is a submodel of $N$ for all $i < \alpha$,
  \item if $i < j < \alpha$ then $N_i \sub N_j$,
  \item $\bigcup_{i < \alpha} N_i = N$.
\end{enumerate}
If, in addition, $\bigcup_{i < \delta} N_i = N_{\delta}$ for every
limit ordinal $\delta < \alpha$, then the sequence is a
\emph{continuous $\alpha$-resolution of $N$}.
\end{defn}

\begin{lem}\label{n:is:ok}
Suppose $\norm{N} = \kappa$ is regular and $\varphi(x; \aaa)$ is
critical for $N$. Then there is an element $c \in \varphi(M; \aaa)$
such that the type $\tp(c/\abs{N}, M)$ is isolated.
\end{lem}
\begin{proof}
Without loss of generality we may assume $\norm{M} = \kappa$. Fix a
club $C = \seq{\alpha_i: i < \kappa} \sub \kappa$ and a continuous
$\kappa$-resolution $\seq{N_i: i < \kappa}$ of $N$ such that
\begin{enumerate}
  \item for all $\alpha_i, \alpha_j \in C$ and $i < j$ we have $\abs{\alpha_i}
  \leq \abs{\alpha_j \setminus \alpha_i}$,
  \item $\norm{N_i} = \abs{\alpha_i}$,
  \item $\aaa \in N_0$.
\end{enumerate}
By the inductive hypothesis we construct a sequence $\seq{b_i \in
\varphi(M; \aaa): i < \kappa}$ such that each type
$\tp(b_i/\abs{N_i}, M)$ is isolated. Fix an enumeration
$\seq{\phi_i: i < \kappa}$ of all the formulas in the language of $T
\cup \ed(N)$ such that for each $\alpha_i \in C$ we have
\[
\set{i: \phi_i \text{ is a formula in the language of } T \cup
\ed(N_i)} \sub \alpha_i.
\]
Now define a function $f: C \fun \kappa$ by letting $f(\alpha_i)$ be
the least ordinal such that $\phi_{f(\alpha_i)}$ isolates the type
$\tp(b_i/\abs{N_i}, M)$. Since $f$ is a pressing-down function on a
stationary subset of $\kappa$ and $\kappa$ is regular, by Fodor's
Lemma, there is a $\gamma < \kappa$ such that $f^{-1}(\gamma) \sub
C$ is stationary. Clearly for any $\alpha_i, \alpha_j \in
f^{-1}(\gamma)$, if $\alpha_i < \alpha_j$ then $\tp(b_i/\abs{N_j},
M) = \tp(b_j/\abs{N_j}, M)$ as they are both isolated by
$\phi_{\gamma}$. So $\tp(b_i/\abs{N}, M) = \tp(b_j/\abs{N}, M)$ for
any $\alpha_i, \alpha_j \in f^{-1}(\gamma)$. And this type is
isolated by $\phi_{\gamma}$ as desired.
\end{proof}

For the case that $\norm{N}$ is singular we need to work harder.
First we formulate the following concept:

\begin{defn}
Let $\seq{N_i: i < \alpha}$ be an $\alpha$-resolution of $N$. Let
$\aaa \in N_0$. Let $\varphi(x; \aaa)$ be critical for $N$. We say
that $\textbf{F} = \seq{\varphi_i(x): i < \alpha}$ is a \emph{spinal
sequence of $\varphi(x; \aaa)$ for $\seq{N_i: i < \alpha}$} if:
\begin{enumerate}
  \item each $\varphi_i(x)$ is a formula in the language of $T \cup \ed(N_i)$,
  \item $\varphi_i(M) \neq \emptyset$ and $\varphi_i(M) \sub \varphi(M; \aaa)$ for each $i <
  \alpha$ ,
  \item if $b \in \varphi_i(M)$ then the type $\tp(b/\abs{N_i}, M)$ is isolated by
  $\varphi_i(x)$.
\end{enumerate}
We write $\dom(\textbf{F})$ for the set
\[
\set{a \in \abs{N}: a \text{ occurs as a parameter in some }
\varphi_i(x) \in \textbf{F}}.
\]
\end{defn}

\begin{lem}\label{countable:cf:ok}
Suppose $\norm{N} = \kappa$ is singular and $\varphi(x; \aaa)$ is
critical for $N$. Then there is an element $c \in \varphi(M; \aaa)$
such that the type $\tp(c/\abs{N}, M)$ is isolated.
\end{lem}
\begin{proof}
As above we may assume $\norm{M} = \kappa$. Let $\lambda =
\cf(\kappa) < \kappa$. Let $\seq{\mu_i: i < \lambda} \sub \kappa$ be
a strictly increasing sequence of cardinals such that it is
unbounded in $\kappa$. Let $\seq{N_i: i < \lambda}$ be a
$\lambda$-resolution of $N$ such that $\aaa \in N_0$ and $\norm{N_i}
= \mu_i$.

Let $\textbf{F}_0$ be a spinal sequence of $\varphi(x; \aaa)$ for
$\seq{N_i: i < \lambda}$. Note that the existence of such a sequence
is guaranteed by the inductive hypothesis. We have
$\abs{\dom(\textbf{F}_0)} \leq \lambda$. Now let $K_0 \sub N$ be the
submodel generated by $\dom(\textbf{F}_0) \cup \set{\aaa}$. Note
that $\varphi(x; \aaa)$ is critical for $K_0$. Since $\norm{K_0}
\leq \lambda < \kappa$, by the inductive hypothesis there is an
element $c_0 \in \varphi(M; \aaa)$ such that $\tp(c_0/\abs{K_0}, M)$
is isolated by some formula $\sigma_0(x)$ in $L(T \cup \ed(K_0))$.
Notice that if $\textbf{F}_0 \sub \tp(c_0/\abs{K_0}, M)$ then we are
done: in this case $\sigma_0(x)$ isolates the entire $\textbf{F}_0$
and each $\varphi_i(x) \in \textbf{F}_0$ isolates the type
$\tp(c_0/\abs{N_i}, M)$, so the type $\tp(c_0/\abs{N}, M)$ is
isolated by $\sigma_0(x)$.

Next, since $\varphi(x; \aaa) \wedge \sigma_0(x)$ is critical for
$N$ (because it contains $\varphi(x; \aaa)$ as a conjunct), we can
find a spinal sequence $\textbf{F}_1$ of $\varphi(x; \aaa) \wedge
\sigma_0(x)$ for $\seq{N_i: i < \lambda}$. Clearly $\textbf{F}_1$ is
also a spinal sequence of $\varphi(x; \aaa)$ for $\seq{N_i: i <
\lambda}$. Let $K_1 \sub N$ be the submodel generated by $\abs{K_0}
\cup \dom(\textbf{F}_1)$. Then, similarly, we can find an element
$c_1 \in \varphi(M; \aaa)$ and a formula $\sigma_1(x)$ in $L(T \cup
\ed(K_1))$ that isolates the type $\tp(c_1/\abs{K_1}, M)$.

Continuing in this fashion we can construct a sequence
$\seq{(\textbf{F}_i, c_i, \sigma_i(x)): i < \lambda^+}$ such that
\begin{enumerate}
  \item $c_i \in \varphi(M; \aaa)$,
  \item $\textbf{F}_{i+1}$ is a spinal sequence of $\varphi(x; \aaa) \wedge
  \sigma_i(x)$ for $\seq{N_i: i < \lambda}$,
  \item $\sigma_i(x)$ is a formula in $L(T \cup
  \ed(K_i))$ which isolates the type $\tp(c_i/\abs{K_i}, M)$, where $K_i \sub N$ is the submodel generated by the
  set $\set{\aaa} \cup \bigcup_{j \leq i} \dom(\textbf{F}_j)$,
  \item if $i$ is a limit ordinal then $\textbf{F}_i$ is not
  defined.
\end{enumerate}
Let $K = \bigcup_{j < \lambda^+} K_j$. Let
\[
S^{\lambda}_{\lambda^+} = \set{\alpha < \lambda^+: \cf(\alpha) =
\lambda},
\]
which is a stationary subset of $\lambda^+$. Fix an enumeration of
all the formulas in $L(T \cup \ed(K))$ such that for each $\alpha
\in S^{\lambda}_{\lambda^+}$ we have
\[
\set{i: \phi_i \text{ is a formula in the language of } T \cup
\ed(K_{\alpha})} \sub \alpha.
\]
So again by Fodor's Lemma there is a $\sigma_j(x)$ and a stationary
subset $S \sub S^{\lambda}_{\lambda^+}$ such that for all $\alpha
\in S$ the type $\tp(c_{\alpha}/\abs{K_{\alpha}}, M)$ is isolated by
$\sigma_j(x)$.

For any $\alpha, \beta \in S$ with $\alpha < \beta$, consider
$\textbf{F}_{\alpha + 1}$. Since $\sigma_{\alpha}(x)$ is
$\sigma_j(x)$, $\textbf{F}_{\alpha+1}$ is a spinal sequence of
$\varphi(x; \aaa) \wedge \sigma_j(x)$ for $\seq{N_i: i < \lambda}$.
So
\[
M \models \ex{x} (\varphi(x; \aaa) \wedge \sigma_j(x) \wedge
\varphi_i(x))
\]
for all $\varphi_i(x) \in \textbf{F}_{\alpha+1}$ (this is by the
second condition in the definition of a spinal sequence above).
Since $\sigma_j(x)$ also isolates the complete type
$\tp(c_{\beta}/\abs{K_{\beta}}, M)$ and $\dom(\textbf{F}_{\alpha+1})
\sub \abs{K_{\beta}}$, we must have $\textbf{F}_{\alpha+1} \sub
\tp(c_{\beta}/\abs{K_{\beta}}, M)$. So $\sigma_j(x)$ isolates
$\textbf{F}_{\alpha+1}$. Since each $\varphi_i(x) \in
\textbf{F}_{\alpha+1}$ determines the type over $N_i$, we see that
$\sigma_j(x)$ isolates the type $\tp(c_{\beta}/\abs{N}, M)$.
\end{proof}

With these two lemmas we can now simply proceed to build an almost
isolating sequence for some model of $T$ over $N$ much in the same
way as in Theorem~\ref{ss:countable:closure}, only now the length of
the almost isolating sequence can go up to $\norm{N} \cdot \omega$.
This proves Theorem~\ref{D:SS:equi}.

We end this paper with a question:

\begin{ques}
Is there an analog of Theorem~\ref{D:SS:equi} for uncountable
languages?
\end{ques}

Notice that, if $T$ is a theory in an uncountable language and the
SS-property and the D-property are not equivalent for $T$, then
there is an $M \models T$ and an $N \sub M$ such that the complete
theory $T \cup \ed(N)$ is not totally transcendental. This is
because primary models always exist for totally transcendental
theories.
%---------------------------------------------------------------------
%Included for Gather Purpose only:
%input "C:\texmf\bibtex\MYbib.bib"
\bibliographystyle{asl}
\bibliography{MYbib}
%---------------------------------------------------------------------

\end{document}